\documentclass[12pt]{amsart} \usepackage{graphicx, psfrag}

\DeclareMathOperator{\Inv}{Inv} 
 
\newcommand{\bs}{\backslash}\newcommand{\ep}{\varepsilon}
 \DeclareMathOperator{\Int}{Int}
\DeclareMathOperator{\cl}{cl}

\theoremstyle{plain} \newtheorem{thm}{Theorem}
 
\newtheorem{lemma}[thm]{Lemma}

\theoremstyle{definition}

\theoremstyle{remark}

\begin{document}
    
\title{positively expansive homeomorphisms of compact spaces} 
\author{David Richeson} \author{Jim Wiseman} \address{Dickinson
College\\
Carlisle, PA 17013} \email{richesod@dickinson.edu} \address{Swarthmore
College\\
Swarthmore, PA 19081} \email{jwisema1@swarthmore.edu}

\date{\today}
\begin{abstract}
We give a new and elementary proof showing that a homeomorphism
$f:X\to X$ of a compact metric space is positively expansive if and only if  $X$ is finite. 
\end{abstract} 

\maketitle
\section{Introduction}\label{sec:intro}

A continuous map $f:X\to X$ on a metric space $X$ is {\em positively
expansive} if there exists $\rho>0$ such that for any distinct $x,y\in
X$ there is an $n\ge 0$ with $d(f^n(x),f^n(y))>\rho$. The constant $\rho$ is called the {\em
expansive constant}.  In this paper we give a simple, new proof of the
following theorem.

\begin{thm}
	\label{thm:fwdexp}
	Let $X$ be a compact metric space.  A homeomorphism $f:X\to X$ is
	positively expansive if and only if $X$ is finite.
\end{thm}

As far as we can ascertain, the first explicit statement of this
theorem was made by Keynes and Robertson (\cite{KR}).  Their proof
used the idea of generators for topological entropy.  Later, the
theorem was proved by Hiraide (\cite{H2}).  His proof requires a
technical result of Reddy's which in turn uses Frink's metrization
theorem to find a compatible metric with respect to which the
homeomorphism is expanding (\cite{R}, \cite{F}, \cite{AH} p.  41).  In
fact, before either of these two papers Gottschalk and Hedlund proved
several results that had, as an unstated corollary, the fact that $X$
must have an isolated point (\cite{GH}, Theorems 10.30, 10.36).  One
can use this observation to prove that all points are isolated, and
thus that $X$ is finite.  In this paper we give a proof that is short
and dynamical and relies only on elementary topological arguments.


As Theorem \ref{thm:fwdexp} illustrates, positive expansiveness is a
very restrictive property.  One cannot restate the theorem for
expansive homeomorphisms (a homeomorphism $f$ is {\em expansive} if
there exists $\rho>0$ such that if $d(f^n(x),f^n(y))<\rho$ for every
integer $n$, then $x=y$).  Although some compact spaces do not admit
expansive homeomorphisms (such as the 2-sphere, the projective plane,
the Klein bottle (\cite{H})), others do.  For instance, O'Brien and
Reddy proved that every compact orientable surface of positive genus
admits an expansive homeomorphism (\cite{OR}).  Also, every Anosov
diffeomorphism is expansive.

Furthermore, one cannot state the same theorem for noninvertible
dynamical systems.  For instance, the doubling map on $S^{1}$ is a
positively expansive continuous map.  Hiraide does prove that no
positively expansive map exists on any manifold with boundary
(\cite{H2}).

We remind the reader of some standard definitions.  Let $f:X\to X$ be
a homeomorphism.  The {\em $\omega$-limit set} of a point $x\in X$ is
defined to be
\[
\omega(x)=\underset{N>0}{\bigcap}\cl\Big(\underset{n>N}{\bigcup}f^n(x)\Big).\]
A set $S$ is {\em invariant} if $f(S)=S$.  We denote the maximal
invariant subset of a set $N$ by $\Inv N$.  An invariant set $S$ is an
{\em isolated invariant set} provided there is a compact neighborhood
$N$ of $S$ with the property that $S=\Inv N$; the set $N$ is an {\em
isolating neighborhood} for $S$.  A set $S$ is an {\em attractor} if
there is an isolating neighborhood $N$ for $S$ with the property that
$f(N)\subset\Int N$; in this case $N$ is called an {\em attracting
neighborhood}.  Likewise, $S$ is a {\em repeller} if it has a {\em
repelling neighborhood}, an isolating neighborhood $N$ with the
property $f^{-1}(N)\subset\Int N$.  Finally, we let $B_{\ep}(x)$ 
denote the $\ep$-ball about $x$.

\section{Bounded dynamical systems}
\label{sec:bounded}

This work relies heavily on the notion of bounded dynamical systems
(see \cite{RW}, \cite{RW2}).  A dynamical system is {\em bounded} if
there exists a compact set $W$ with the property that the forward
orbit of every point in $X$ intersects $W$.  Such a set, $W$, is
called a {\em window}.  Clearly every dynamical system on a compact
space $X$ is bounded, thus the notion of boundedness is only
interesting for noncompact spaces.

Below we state several properties that are equivalent to boundedness;
the theorem is proved in \cite{RW}, but since the proof is short we
include it again here.  We note that the theorem is also true for
flows or semiflows and the proof is nearly identical to the one given
below.

\begin{thm}\label{thm:bounded} If $X$ is a locally compact metric
space and $f:X\to X$ is a continuous map, then the following are
equivalent.
\begin{enumerate}
	\item
	$f$ is bounded.
	\item
	There is a compact set $V$ such that $\emptyset\ne\omega(x)\subset
	V$ for all $x\in X$.
	\item
	There exists a forward invariant window.
	\item
	There is a compact global attractor $\Lambda$ (that is, there is
	an attractor $\Lambda$ with the property that
	$\emptyset\ne\omega(x)\subset\Lambda$ for every $x\in X$).
\end{enumerate}
\end{thm}
\begin{proof}
	It is clear that (4)$\implies$(3)$\implies$(2)$\implies$(1). 
	Thus, we must prove that the existence of a window implies the
	existence of a compact global attractor, (1)$\implies$(4).

	Suppose $f$ has a window $W$.  It suffices to show that there is a
	window $W_{1}$ with the property $f(W_{1})\subset \Int(W_{1})$. 
	For each $x\in X$ there is an $n_{x}\ge 0$ such that
	$f^{n_{x}}(x)\in W$.  Let $\delta>0$, and let
	$W_{0}=\cl(B_{\delta}(W))$, the closure of the
	$\delta$-neighborhood of $W$.  Clearly, for each $x\in W_{0}$
	$\cl(B_{\delta/2}(f^{n_x}(x)))\subset\Int W_{0}$.  Moreover, there
	is an open neighborhood $U_x$ of $x$ such that
	$\cl(B_{\delta/2}(f^{n_x}(y)))\subset\Int W_{0}$ for all $y\in
	U_x$.  The sets $\{U_x:x\in W_{0}\}$ form an open cover of
	$W_{0}$.  Since $W_{0}$ is compact there is a finite subcover,
	$\{U_{x_1},\ldots,U_{x_m}\}$.  Let
	$n=\max\{n_{x_k}:k=1,\ldots,m\}$.  It follows that
	$\bigcup_{k=0}^{n}f^k(W_{0})$ is a forward invariant window (thus
	proving (3)).  However, we would like the stronger result of (4).

	Consider the multivalued map $V_{r}(x)=B_{r}(x)$.  By the
	compactness of $W_{0}$, there is an $\ep>0$ such that
	$(V_{\ep}\circ f)^{n_{x_{i}}}(y)\subset\Int W_{0}$ for all $y\in
	U_{x_{i}}$.  Then, the set $W_{1}=\bigcup_{k=0}^{n}(V_{\ep}\circ
	f)^k(W_{0})$ has the desired property.
\end{proof}

\section{positively expansive homeomorphisms on compact spaces}
\label{sec:compact}

In the discussion that follows it is necessary to work in the product
space $X\times X$.  Given a homeomorphism $f:X\to X$ we use the
notation $f\times f:X\times X\to X\times X$ to denote the
homeomorphism $(f\times f)(x_{1},x_{2})=(f(x_{1}),f(x_{2}))$.  Also,
we let $\Delta=\{(x,x): x\in X\}$ denote the {\em diagonal} of
$X\times X$.

It is well known that a homeomorphism $f:X\to X$ of a compact space
$X$ is expansive if and only if the diagonal $\Delta$ is an isolated
invariant set for $f\times f$ (\cite{A}). 
Analogously we prove that for positively expansive homeomorphisms the
diagonal is a repeller. 

\begin{lemma}\label{lem:diagattractor}
	Let $f:X\to X$ be a positively expansive homeomorphism 
	of a
	compact space $X$.  Then $\Delta$ is a repeller for $f\times
	f:X\times X\to X\times X$.
\end{lemma}
\begin{proof}
	Suppose $X$ is a compact space and $f:X\to X$ is a positively
	expansive homeomorphism with expansive constant $\rho$.  If $X$ is
	a one-point space, the conclusion of the lemma is clearly true. 
	Thus we may assume that $X$ has at least two points.  Consider the
	space $X\times X$ and the homeomorphism $F=f\times
	f$. $F$ restricts
	to a homeomorphism $F_{Y}:Y\to Y$ where $Y=(X\times X)\bs\Delta$. 
	Let $W=\{(x,y)\in Y: d_{X}(x,y)\ge \rho\}$.  Clearly $W$ is a
	compact set and, since $f$ is positively expansive, the forward
	orbit of every point in $Y$ intersects $W$.  Thus $W$ is a window
	for $F_{Y}$, and we conclude that $F_{Y}$ is bounded.

	By Theorem \ref{thm:bounded} there exists a window $W_1\subset Y$
	for $F_{Y}$ with $F_{Y}(W_{1})\subset\Int (W_{1})$.  Then the set
	$N=\cl((X\times X)\bs W_1)$ has the property that
	$F^{-1}(N)\subset\Int N$ and $\Inv N=\Delta$.  Thus $\Delta$ is
	a repeller for $F$.
\end{proof}
  
\begin{proof}[Proof of Theorem \ref{thm:fwdexp}]
	Let $f:X\to X$ be a positively expansive homeomorphism of a
	compact space $X$.  Let $g=f^{-1}$ and $G=g\times g:X\times X\to
	X\times X$.  By Lemma \ref{lem:diagattractor} the diagonal
	$\Delta\subset X\times X$ is an attractor for $G$.  Thus, for
	$(x,y)$ sufficiently close to $\Delta$, $G^n(x,y)\to\Delta$ as
	$n\to\infty$.  More precisely, there exists $\ep>0$ such that if
	$d(x,y)<\ep$ then $d(g^n(x),g^n(y))\to 0$ as $n\to\infty$.

	Define an equivalence relation $\sim$ on $X$ as follows: $x\sim y$
	iff there exists a sequence of points $x=x_0,x_1,...,x_r=y$ such
	that $d(x_k,x_{k+1})<\ep$ for $k=0,...,r-1$.  This equivalence
	relation is an open condition, thus each equivalence class is an
	open subset of $X$.  Since the set of equivalence classes is a
	cover of $X$ by mutually disjoint open sets, the compactness of
	$X$ implies that there are only finitely many, $U_1,...,U_m$. 
	Also, since each $U_i$ contains its limit points, it is closed,
	and hence compact.

	Let $U$ be an equivalence class, and let $x,y\in U$.  Then there
	is a sequence $x=x_0,x_1,...,x_r=y$ such that $d(x_k,x_{k+1})<\ep$
	for $k=0,...,r-1$.  So, $d(g^{n}(x_{k}),g^{n}(x_{k+1}))\to 0$ as
	$n\to\infty$ for each $k=0,...,r-1$.  Thus
	$d(g^{n}(x),g^{n}(y))\to 0$ as $n\to\infty$.  Since $U$ is
	compact, the diameter of $g^{n}(U)$ goes to zero as $n\to\infty$.

	For each $n$, $g^{n}(U_{1}),...,g^{n}(U_{m})$ is a collection of
	mutually disjoint sets whose union is all of $X$.  Moreover, the
	diameter of each set $g^{n}(U_{i})$ can be made arbitrarily small
	(letting $n$ get large).  Thus, it must be the case that each
	$U_{i}$ consists of a single point, and that $X$ is a finite set. 
\end{proof}

\bibliographystyle{amsplain} 
\bibliography{posexpcompact} 
\nocite{A}
\nocite{C}
\end{document}